O.M. Mokliachuk*

Igor Sikorsky Kyiv Polytechnic Institute, Kyiv, Ukraine

# ESTIMATION OF ACCURACY AND RELIABILITY OF MODELS OF φ-SUB-GAUSSIAN STOCHASTIC PROCESSES IN $C(T)$ SPACES

**Background.** At present, in the theory of stochastic process modeling a problem of assessment of reliability and accuracy of stochastic process model in $C(T)$ space wasn't studied for the case of inexplicit decomposition of process in the form of a series with independent terms.
**Objective.** The goal is to study reliability and accuracy in $C(T)$ of models of processes from $\text{Sub}_\varphi(\Omega)$ that cannot be decomposed in a series with independent elements explicitly.
**Methods.** Using previous research in the field of modeling of stochastic processes, assumption is considered about possibility of decomposition of a stochastic process in the series with independent elements that can be found using approximations.
**Results.** Impact of approximation error of process decomposition in series with independent elements on reliability and accuracy of modeling of stochastic process in $C(T)$ is studied.
**Conclusions.** Theorems are proved that allow estimation of reliability and accuracy of a model in $C(T)$ of a stochastic process from $\text{Sub}_\varphi(\Omega)$ in the case when decomposition of this process in a series with independent elements can be found only with some error, for example, using numerical approximations.
**Keywords:** stochastic processes; φ-sub-Gaussian processes; models of stochastic processes; reliability and accuracy of models of stochastic processes.

## Introduction

In the theory of stochastic processes, modeling of stochastic processes is a vital part. This problem has been a matter of active research in recent decades. It has become an irreplaceable part of research, development and applications in many different fields. As a result, to build the model of a stochastic process and study its properties is the task of high importance.

Much attention has been paid to the models of stochastic processes from $\text{Sub}_\varphi(\Omega)$ spaces, particularly in [1, 2]. One of the most used spaces for model development of stochastic processes is $C(T)$ space, the space of functions, continuous on the interval $T$. Models of stochastic processes in these spaces are considered, for example, in [3] and [4]. The following representation of this stochastic process is used [4] to build a stochastic process model with given reliability and accuracy in any space:

$$X(t) = \sum_{k=0}^{\infty} \xi_k a_k(t).$$

Here, $\xi_k$ are random variables, $a_k(t)$ are functions, and this series has to be mean square convergent in the space where we want to build the model of this process.

The model of the process itself is defined as the sum of first $N$ elements of the above series and is usually denoted by $X_N(t)$ [4], i.e.

$$X_N(t) = \sum_{k=0}^{N} \xi_k a_k(t). \qquad (1)$$

One of the biggest problems of dealing with such decompositions is that it is often impossible to find functions $a_k(t)$ in an explicit way, and their approximations have to be used. Needless to say, errors of these approximations will have the influence on the reliability and accuracy of the developed model. This problem has been previously considered in [5–8].

## Research objective

In this paper, we consider models of stochastic processes in $C(T)$ space when decomposition elements $a_k(t)$ cannot be found explicitly and prove theorems that allow modeling of such stochastic process taking into account errors of approximations of $a_k(t)$.

## Basic concept

Let $(\Omega, \mathcal{F}, P)$ be a standard probability space, $T$ be some parametric space.

---

* corresponding author: omoklyachuk@gmail.com

***Definition* 1** [9]. Function $\varphi(x)$ is called an Orlicz N-function, if it is even convex continuous function and

$$\varphi(0) = 0, \quad \varphi(x) > 0, \quad x \neq 0;$$
$$\varphi(x)/x \to 0, \quad x \to 0;$$
$$\varphi(x)/x \to \infty, \quad x \to \infty.$$

**Assumption Q** [10]. Let $\varphi(x)$ be an Orlicz N-function and $\lim_{x \to 0} \frac{\varphi(x)}{x^2} = c > 0$.

The constant $c$ can be equal to $+\infty$.

***Definition* 2.** Random variable $\xi$ is called $\varphi$-sub-Gaussian, if there exists such $a \geq 0$, that for all $\lambda \in R$ inequality

$$E \exp\{\lambda \xi\} \leq \exp\{\varphi(a\lambda)\},$$

holds true, where $\varphi$ is N-function that satisfies assumption Q.

Function $\tau(\xi)$ defined as

$$\tau(\xi) = \inf\{a \geq 0 : E \exp\{\lambda \xi\} \leq \exp\{\varphi(a\lambda)\}, \lambda \in R\}$$

is called the $\varphi$-sub-Gaussian standard of $\xi$.

A space of $\varphi$-sub-Gaussian random variables is usually denoted as $\text{Sub}_\varphi(\Omega)$. This space is Banach space with the norm $\tau(\xi)$ and $E\xi = 0$.

***Definition* 3.** Stochastic process $X = \{X(t), t \in T\}$ is called $\varphi$-sub-Gaussian, if random variables $X(t)$ for each $t \in T$ are $\varphi$-sub-Gaussian.

The next statement will be required for further proofs.

**Theorem 1** [10]. Let $\xi_1, \xi_1, \ldots, \xi_n \in \text{Sub}_\varphi(\Omega)$ be independent random variables. If function $\varphi(|x|^{1/s})$, $x \in R$, is convex for $s \in (0, 2]$, then

$$\tau_\varphi^s \left( \sum_{k=1}^n \xi_k \right) \leq \sum_{k=1}^n \tau_\varphi^s(\xi_k).$$

Let us make some assumptions. Let $X = \{X(t), t \in B\}$ be a continuous with probability 1 stochastic process from $\text{Sub}_\varphi(\Omega)$ space (conditions for continuity can be found in [1]), let $(B, \rho)$ be a compact, and let $X$ be separable in $(B, \rho)$. Assume that continuous monotone increasing function $\sigma = \{\sigma(h), h > 0\}$ exists, $\sigma(h) \to 0$ when $h \to 0$, and next inequality holds true:

$$\sup_{\rho(t,s) \leq h} \tau_\varphi(X(t) - X(s)) \leq \sigma(h). \qquad (2)$$

Besides, let $\beta > 0$ be the following number:

$$\beta = \sigma\left(\inf_{s \in B} \sup_{t \in B} \rho(t, s)\right).$$

***Definition* 4.** Minimal number of closed spheres with the radius $u$ that can cover $(B, \rho)$ is called metric massiveness $N_B(u)$.

Under these assumptions, next modification of the theorem from [1, p. 61] holds true.

**Theorem 2.** Let stochastic process $X = \{X(t), t \in B\}$ belong $\text{Sub}_\varphi(\Omega)$ space, $(B, \rho)$ — be the compact, and let $X$ be separable in $(B, \rho)$. Assume that condition (2) holds true for the process $X$, $r = \{r(u) : u \geq 1\}$ is the continuous function such that $r(u) > 0$ for $u > 1$, and function $s(t) = r(\exp\{t\})$, $t > 0$, is convex. Then, if

$$\int_0^\beta r(N_B(\sigma^{(-1)}(u)))du < \infty,$$

stochastic process $X(t)$ is bounded with probability 1, and for all $p \in (0, 1)$ and $x > 0$ next inequalities are true:

$$P\{\sup_{t \in (0,T)} X(t) > x\} \leq Z_r(p, \beta, x),$$
$$P\{\inf_{t \in (0,T)} X(t) < -x\} \leq Z_r(p, \beta, x),$$
$$P\{\sup_{t \in (0,T)} |X(t)| > x\} \leq 2Z_r(p, \beta, x),$$

$$Z_r(p, \beta, x) = \inf_{\lambda > 0} \exp\left\{\theta_\varphi(\lambda, p) + p\varphi\left(\frac{\lambda \beta}{1-p}\right) - \lambda x\right\}$$
$$\times r^{(-1)}\left(\frac{1}{\beta p} \int_0^{\beta p} r(N_B(\sigma^{(-1)}(u)))du\right),$$

where

$$\theta_\varphi(\lambda, p) = \sup_{u \in B}\left((1-p)\varphi\left(\frac{\gamma(u)\lambda}{1-p}\right)\right),$$
$$\gamma(u) = \tau_\varphi(X(u)).$$

This theorem is very general, so let us consider cases where function $\varphi(t)$ is given particularly.

**Theorem 3.** Let stochastic process Let $X = \{X(t), t \in B\}$ belong to $\text{Sub}_\varphi(\Omega)$, $(B, \rho)$ — be the compact and let $X$ be separable on $(B, \rho)$. Let $\varphi(t) = \frac{t^\zeta}{\zeta}$, $\zeta \geq 2$, $t > 1$, and let process $X$ satisfy the conditions of the previous theorem. Then the process $X(t)$ is bounded with probability 1, and for all $p \in (0,1)$ next inequalities hold true:

$$P\{\sup_{t \in (0,T)} X(t) > x\} \leq Z_{r_1}(p, \beta, x),$$

$$P\{\inf_{t \in (0,T)} X(t) < -x\} \leq Z_{r_1}(p, \beta, x),$$

$$P\{\sup_{t \in (0,T)} |X(t)| > x\} \leq 2Z_{r_1}(p, \beta, x),$$

$$Z_{r_1}(p, \beta, x) = \exp\left\{\frac{(1-\zeta)(x(1-p))^{\frac{\zeta}{\zeta-1}}}{\zeta(\gamma^\zeta(1-p) + p\beta^\zeta)^{\frac{1}{\zeta-1}}}\right\}$$

$$\times r_1^{(-1)}\left(\frac{1}{\beta p}\int_0^{\beta p} r_1(N_B(\sigma^{(-1)}(u)))du\right),$$

$$x > \frac{\gamma^\zeta(1-p) + p\beta^\zeta}{(1-p)v^{\zeta-1}},$$

where $v = \min\{\beta, \gamma\}$, $\gamma^\zeta = \sup_{u \in B}\tau_\varphi^\zeta(X(u))$.

P r o o f. Let

$$(1-p) \leq \gamma\lambda \quad \text{and} \quad (1-p) \leq \beta\lambda. \quad (3)$$

Consider the power of the exponent in the formula of $Z_{r_1}(p, x, \beta)$ from the previous theorem using $\varphi(t)$ and $\theta_\varphi(\lambda, p)$ provided in the current theorem. Then

$$\sup_{u \in (0,T)}\left((1-p)\varphi\left(\frac{\gamma(u)\lambda}{1-p}\right)\right) + p\varphi\left(\frac{\lambda\beta}{1-p}\right) - \lambda x$$

$$= \sup_{u \in B}\left(\frac{\gamma^\zeta(u)\lambda^\zeta}{\zeta(1-p)^{\zeta-1}}\right) + p\frac{\lambda^\zeta\beta^\zeta}{\zeta(1-p)^\zeta} - \lambda x$$

$$= \frac{\lambda^\zeta}{\zeta}\left(\frac{\sup_{u \in B}\gamma^\zeta(u)(1-p) + p\beta^\zeta}{(1-p)^\zeta}\right) - \lambda x. \quad (4)$$

If $\gamma^\zeta := \sup_{u \in B}\gamma^\zeta(u) = \sup_{u \in B}\tau_\varphi^\zeta(X(u))$, formula (4) reaches its minimum when

$$\lambda = \left(\frac{x(1-p)^\zeta}{\gamma^\zeta(1-p) + p\beta^\zeta}\right)^{\frac{1}{\zeta-1}},$$

and minimal value of formula (4) is equal to

$$\frac{(1-\zeta)(x(1-p))^{\frac{\zeta}{\zeta-1}}}{\zeta(\gamma^\zeta(1-p) + p\beta^\zeta)^{\frac{1}{\zeta-1}}}.$$

Besides, if $\lambda$ minimizes (4), inequalities (3) hold true, or

$$x > \frac{\gamma^\zeta(1-p) + p\beta^\zeta}{(1-p)v^{\zeta-1}},$$

where $v = \min\{\beta, \gamma\}$. □

**Models of stochastic processes**

Let us define model of the stochastic process as it is given in (1).

***Definition 5.*** Stochastic process $X_N(t) = \{X_N(t), t \in T\}$ is called the model of stochastic process $X = \{X(t), t \in T\}$, if

$$X(t) = \sum_{k=0}^{\infty} \xi_k a_k(t)$$

and

$$X_N(t) = \sum_{k=0}^{N} \xi_k a_k(t).$$

***Definition 6.*** The model $X_N(t)$ approximates the process $X(t)$ with given reliability $1 - v$ and accuracy $\delta$ in $C(B)$ space, if

$$P\{\sup_{t \in T}|\Delta_N(t)| > \delta\} \leq v,$$

where

$$\Delta_N(t) = X(t) - X_N(t).$$

Assume that continuous monotone increasing function $\sigma_N = \{\sigma_N(h), h > 0\}$ exists, $\sigma_N(h) \to 0$, $h \to 0$, and let the next inequality holds true:

$$\sup_{\rho(t,s) \leq h} \tau_\varphi(X_N(t) - X_N(s)) \leq \sigma_N(h). \quad (5)$$

Next two statements allow estimation of reliability and accuracy of a model of the stochastic processes in $C(T)$ space for two different classes of function $\varphi(t)$.

Let us consider two classes of functions. Let denote two classes of functions [9]:

$$C_{\mathrm{I}} = \left\{\varphi(t) \mid \varphi(t) = \frac{t^\zeta}{\zeta}, \zeta \geq 2, t > 1\right\},$$

$$C_{\mathrm{II}} = \left\{\varphi(t) \mid \varphi(t) = \frac{t^\zeta}{\zeta}, 1 < \zeta \leq 2\right\}.$$

**Theorem 4.** Let stochastic process $X_N(t) = \{X_N(t), t \in T\}$ belong to $\mathrm{Sub}_\varphi(\Omega)$ space, $(B, \rho)$ — be the compact and let $X$ be separable on $(B,\rho)$. Let $\varphi(t) \in C_I$, and let the process $X$ satisfy the conditions of theorem 2. Then, the model $X_N(t)$ approximates process $X(t)$ with the given reliability $1-\nu$ and accuracy $\delta$ in $C(B)$ space, if

$$\nu \le 2\exp\left\{-\frac{(\zeta-1)(\delta(1-p))^{\frac{\zeta}{\zeta-1}}}{\zeta(\gamma_N^\zeta(1-p)+p\beta^\zeta)^{\frac{1}{\zeta-1}}}\right\}r_1^{(-1)}$$

$$\times\left(\frac{1}{\beta p}\int_0^{\beta p}r_1(N_B(\sigma_N^{(-1)}(u)))du\right),$$

$$\delta > \frac{\gamma_N^\zeta(1-p)+p\beta^\zeta}{(1-p)v^{\zeta-1}},$$

where $v = \min\{\beta, \gamma_N\}$, $\gamma_N^\zeta = \sup_{u\in[0,T]}\tau_\varphi^\zeta(\Delta_N(u))$.

P r o o f. This theorem is the direct corollary of theorem 3. □

**Theorem 5.** Let stochastic process $X_N(t) = \{X_N(t), t \in T\}$ belong $\text{Sub}_\varphi(\Omega)$ space, let $(B,\rho)$ be the compact and let $X$ be separable on $(B,\rho)$. Let $\varphi(t) = \frac{t^\zeta}{\zeta}$, $1 < \zeta \le 2$, and let the process $X$ satisfy the conditions of the theorem 2. Then, the model $X_N(t)$ approximates the process $X(t)$ with given reliability $1-\nu$ and accuracy $\delta$ in $C(B)$ space, if

$$\nu \le 2\exp\left\{-\frac{(\zeta-1)(\delta(1-p))^{\frac{\zeta}{\zeta-1}}}{\zeta(\gamma_N^\zeta(1-p)+p\beta^\zeta)^{\frac{1}{\zeta-1}}}\right\}r_1^{(-1)}$$

$$\times\left(\frac{1}{\beta p}\int_0^{\beta p}r_1(N_B(\sigma_N^{(-1)}(u)))du\right),$$

where $\gamma_N^\zeta = \sup_{u\in[0,T]}\tau_\varphi^\zeta(\Delta_N(u))$.

P r o o f. Current theorem is the corollary of theorem 3. In this case we do not need any restrictions on $\delta$. □

Let us now consider the case, when we have some specific function $\sigma(h)$ for condition (2), namely $\sigma(h) = Ch^æ$, $C$ — some non-negative constant, $æ \in (0,1)$.

**Theorem 6.** Let stochastic process $X = \{X(t), t \in [0,T]\}$ belong to $\text{Sub}_\varphi(\Omega)$ space, $\varphi(t) = \frac{t^\zeta}{\zeta}$, $\zeta \ge 2$, $t > 1$. Let $X$ be separable, and let it satisfy condition (2), where $\sigma(h) = Ch^æ$, $æ \in (0,1)$. Then, process $X(t)$ is bounded with probability 1 and

$$P\{\sup_{t\in(0,T)} X(t) > x\} \le Z_{r_1}(x),$$

$$P\{\inf_{t\in(0,T)} X(t) < -x\} \le Z_{r_1}(x),$$

$$P\{\sup_{t\in(0,T)} |X(t)| > x\} \le 2Z_{r_1}(x),$$

where

$$Z_{r_1}(x)$$

$$= \exp\left\{-(x-\gamma)^{\frac{\zeta}{\zeta-1}}\frac{(\zeta-1)x^{\frac{1}{\zeta-1}}}{\zeta(\gamma^\zeta(x-\gamma)+\beta^\zeta\gamma)^{\frac{1}{\zeta-1}}}\right\}2(ex)^{1/æ},$$

$$x > \frac{\gamma(v^{\zeta-1}+1)+\sqrt{\gamma^2(v^{\zeta-1}+1)^2 + 4v^{\zeta-1}(C^\zeta(T/2)^{æ\zeta}-\gamma^2)}}{2v^{\zeta-1}},$$

when $v = \min\{C(T/2)^æ, \gamma\}$, $\gamma^\zeta = \sup_{u\in[0,T]}\tau_\varphi^\zeta(X(u))$.

P r o o f. This theorem follows from theorem 3. We can take function $r_1(t) = t^\alpha$, $0 < \alpha < æ$ as $r_1$. Under the conditions of current theorem, the second part of the function $Z_{r_1}(p,\beta,x)$ from theorem 3 can be transformed in the next way:

$$r_1^{(-1)}\left(\frac{1}{\beta p}\int_0^{\beta p}r_1(N_B(\sigma^{(-1)}(u)))du\right)$$

$$\le r_1^{(-1)}\left(\frac{1}{\beta p}\int_0^{\beta p}r_1\left(\frac{T}{2\sigma^{(-1)}(u)}+1\right)du\right).$$

Besides, metric massivity satisfies the inequality $N_B(u) \le \frac{T}{2u}+1$ on $[0,T]$. Because of $u \le \beta p \le \beta = \sigma(T/2)$, we have $\sigma^{(-1)}(u) \le T/2$ for such $u$. As the result, $T/2\sigma^{(-1)} \ge 1$, and

$$r_1^{(-1)}\left(\frac{1}{\beta p}\int_0^{\beta p}r_1\left(\frac{T}{2\sigma^{(-1)}(u)}+1\right)du\right)$$

$$\le \left(\frac{1}{\beta p}\int_0^{\beta p}\left(\frac{T}{\sigma^{(-1)}(u)}\right)^\alpha du\right)^{1/\alpha}.$$

Using the formula for $\sigma(u)$, we obtain

$$\left(\frac{1}{\beta p}\int_0^{\beta p}\left(\frac{T}{\sigma^{(-1)}(u)}\right)^\alpha du\right)^{1/\alpha} = \left(\frac{1}{\beta p}\int_0^{\beta p}\frac{T^\alpha C^{\alpha/\ae}}{u^{\alpha/\ae}}du\right)^{1/\alpha}$$

$$= \left(\frac{1}{\beta p}\frac{T^\alpha C^{\alpha/\ae}(\beta p)^{1-\frac{\alpha}{\ae}}}{1-\alpha/\ae}\right)^{1/\alpha} \leq \frac{TC^{1/\ae}}{(\beta p)^{1/\ae}(1-\alpha/\ae)^{1/\alpha}}.$$

When $\alpha \to 0$ we will have

$$\frac{TC^{1/\ae}}{(\beta p)^{1/\ae}(1-\alpha/\ae)^{1/\alpha}} \to \frac{TC^{1/\ae}}{(\beta p)^{1/\ae}}e^{1/\ae}.$$

Because of $\beta = \sigma(\inf_{t\in(0,T)}\sup_{t\in(0,T)}\rho(t,s))$ $= C(T/2)^\ae$, we have

$$\frac{TC^{1/\ae}}{(\beta p)^{1/\ae}}e^{1/\ae} = 2\left(\frac{e}{p}\right)^{1/\ae}.$$

As the result,

$$Z_{r_1}(p,\beta,x) = \exp\left\{\frac{(1-\zeta)(x(1-p))^{\frac{\zeta}{\zeta-1}}}{\zeta(\gamma^\zeta(1-p)+p\beta^\zeta)^{\frac{1}{\zeta-1}}}\right\}2\left(\frac{e}{p}\right)^{1/\ae}.$$

The statement of the theorem is derived from the last equality if we denote $p = \gamma/x$. Moreover, condition for $x$ of the previous theorem takes the form

$$x > \frac{\gamma(1-1/x)+(1/x)C^\zeta(T/2)^{\ae\zeta}}{(1-1/x)v^{\zeta-1}},$$

therefore

$$x > \frac{\gamma(v^{\zeta-1}+1)+\sqrt{\gamma^2(v^{\zeta-1}+1)^2+4v^{\zeta-1}(C^\zeta(T/2)^{\ae\zeta}-\gamma^2)}}{2v^{\zeta-1}},$$

where $v = \min\{\beta,\gamma\}$. □

Following corollaries provide reliability and accuracy for models of stochastic processes in $\mathrm{Sub}_\varphi(\Omega)$, when $\varphi$ belongs to classes $C_\mathrm{I}$ or $C_\mathrm{II}$.

**Theorem 7.** Let stochastic process $X = \{X(t), t \in [0,T]\}$ belong to $\mathrm{Sub}_\varphi(\Omega)$ space, $\varphi(t) \in C_\mathrm{I}$. Let $X$ be separable and satisfy condition (2), where $\sigma(h) = Ch^\ae$. Then, the model $X_N(t)$ approximates the process $X(t)$ with given reliability $1-v$ and accuracy $\delta$ in $C(0,T)$ space, if

$$v \leq 2\exp\left\{-\frac{(\zeta-1)\delta^{\frac{1}{\zeta-1}}(\delta-\gamma_N)^{\frac{\zeta}{\zeta-1}}}{\zeta(\gamma_N^\zeta(\delta-\gamma_N)+\beta^\zeta\gamma_N)^{\frac{1}{\zeta-1}}}\right\}2(e\delta)^{1/\ae},$$

when

$$\delta > \frac{\gamma_N(v^{\zeta-1}+1)+\sqrt{\gamma_N^2(v^{\zeta-1}+1)^2+4v^{\zeta-1}(C^\zeta(T/2)^{\ae\zeta}-\gamma_N^2)}}{2v^{\zeta-1}},$$

$v = \min\{\beta,\gamma_N\}$, $\gamma_N^\zeta = \sup_{u\in[0,T]}\tau_\varphi^\zeta(\Delta_N(u))$.

P r o o f. This theorem is the direct corollary of theorem 6. □

**Theorem 8.** Let stochastic process $X = \{X(t), t \in [0,T]\}$ belong to $\mathrm{Sub}_\varphi(\Omega)$ space, $\varphi(t) \in C_\mathrm{II}$. Let $X$ be separable and let it satisfy condition (2), where $\sigma(h) = Ch^\ae$. Then, the model $X_N(t)$ approximates the process $X(t)$ with given reliability $1-v$ and accuracy $\delta$ in $C(0,T)$ space, if

$$v \leq 2\exp\left\{-\frac{(\zeta-1)\delta^{\frac{1}{\zeta-1}}(\delta-\gamma_N)^{\frac{\zeta}{\zeta-1}}}{\zeta(\gamma_N^\zeta(\delta-\gamma_N)+\beta^\zeta\gamma_N)^{\frac{1}{\zeta-1}}}\right\}2(e\delta)^{1/\ae},$$

$$\gamma_N^\zeta = \sup_{u\in[0,T]}\tau_\varphi^\zeta(\Delta_N(u)).$$

P r o o f. This statement is the corollary of theorem 6. Here no restrictions on $\delta$ are needed. □

**Models of stochastic processes that allow representation in series with independent elements**

Assume that stochastic process $X_N(t) = \{X_N(t), t \in T\}$ can be represented as

$$X(t) = \sum_{k=1}^\infty a_k(t)\xi_k, \qquad (6)$$

where $\xi_k \in \mathrm{Sub}_\varphi(\Omega)$. Let $\delta_k(t) = |a_k(t)-\hat{a}_k(t)|$, $\hat{a}_k(t)$ be approximations of $a_k(t)$, and let $\sigma_N$ be the same as in definition 5. Let $X_N$, sum of the first $N$ elements of this series, be the model of such process, as it is provided in definition 5.

For given stochastic process, we can prove next theorems.

**Theorem 9.** Let stochastic process $X_N(t) = \{X_N(t), t \in T\}$ belong to $\mathrm{Sub}_\varphi(\Omega)$ space, $\varphi(t) \in C_\mathrm{I}$, let $(B,\rho)$ be a compact and let $X$ be separable in $(B,\rho)$. Assume that process $X$ satisfy conditions of theorem 2.

Then, the model $X_N(t)$ approximates the process $X(t)$ with given reliability $1 - v$ and accuracy $\delta$ in $C(B)$ space, if

$$v \leq 2\exp\left\{-\frac{(\zeta-1)(\delta(1-p))^{\frac{\zeta}{\zeta-1}}}{\zeta(\gamma_N^\zeta(1-p)+p\beta^\zeta)^{\frac{1}{\zeta-1}}}\right\}r_1^{(-1)}$$
$$\times\left(\frac{1}{\beta p}\int_0^{\beta p}r_1(N_B(\sigma_N^{(-1)}(u)))du\right),$$
$$\delta > \frac{\gamma_N^\zeta(1-p)+p\beta^\zeta}{(1-p)v^{\zeta-1}},$$

where $v = \min\{\beta, \gamma_N\}$, $\gamma_N^\zeta = \sup_{u\in B}\tau_\varphi^\zeta(\Delta_N(u))$,

$$\gamma_N^\zeta \leq \left(\sum_{k=1}^N \tau_\varphi^2(\xi_k)\sup_{u\in B}\delta_k^2(u) + \sum_{k=N+1}^\infty \tau_\varphi^2(\xi_k)\sup_{u\in B}a_k^2(u)\right)^{\zeta/2}.$$

P r o o f.

$$\gamma_N^\zeta = \sup_{u\in B}\tau_\varphi^\zeta(\Delta_N(u))$$
$$= \sup_{u\in B}\tau_\varphi^\zeta\left(\sum_{k=1}^N \xi_k\delta_k(u) + \sum_{k=N+1}^\infty \xi_k a_k(u)\right)$$
$$\leq \left(\sum_{k=1}^N \tau_\varphi^2(\xi_k)\sup_{u\in B}\delta_k^2(u) + \sum_{k=N+1}^\infty \tau_\varphi^2(\xi_k)\sup_{u\in B}a_k^2(u)\right)^{\zeta/2}.$$

Last inequality follows from the properties of function $\tau_\varphi$ and theorem 1. □

**Theorem 10.** Let stochastic process $X_N(t) = \{X_N(t), t \in T\}$ belong to $\mathrm{Sub}_\varphi(\Omega)$ space, $\varphi(t) \in C_{II}$ let $(B, \rho)$ be a compact and let $X$ be separable in $(B, \rho)$. Assume that process $X$ satisfies conditions of theorem 2. Then, the model $X_N(t)$ approximates the process $X(t)$ with given reliability $1 - v$ and accuracy $\delta$ in $C(B)$ space, if

$$v \leq 2\exp\left\{-\frac{(\zeta-1)(\delta(1-p))^{\frac{\zeta}{\zeta-1}}}{\zeta(\gamma_N^\zeta(1-p)+p\beta^\zeta)^{\frac{1}{\zeta-1}}}\right\}r_1^{(-1)}$$
$$\times\left(\frac{1}{\beta p}\int_0^{\beta p}r_1(N_B(\sigma_N^{(-1)}(u)))du\right),$$

where $v = \min\{\beta, \gamma_N\}$,

$$\gamma_N^\zeta \leq \sum_{k=1}^N \tau_\varphi^\zeta(\xi_k)\sup_{u\in B}\delta_k^\zeta(u) + \sum_{k=N+1}^\infty \tau_\varphi^\zeta(\xi_k)\sup_{u\in B}a_k^\zeta(u).$$

P r o o f.

$$\gamma_N^\zeta = \sup_{u\in B}\tau_\varphi^\zeta(\Delta_N(u))$$
$$= \sup_{u\in B}\tau_\varphi^\zeta\left(\sum_{k=1}^N \xi_k\delta_k(u) + \sum_{k=N+1}^\infty \xi_k a_k(u)\right)$$
$$\leq \sum_{k=1}^N \tau_\varphi^\zeta(\xi_k)\sup_{u\in B}\delta_k^\zeta(u) + \sum_{k=N+1}^\infty \tau_\varphi^\zeta(\xi_k)\sup_{u\in B}a_k^\zeta(u).$$

The last inequality follows from the properties of function $\tau_\varphi$ and theorem 1. □

Let the process $X$ be defined on $[0,T]$, and let it satisfy the next condition:

**(C1)** $\sigma(h) = Ch^\ae$, $\delta_k(t) - \delta_k(s) \leq \widehat{C}_k h^\ae$, $a_k(t) - a_k(s)a_k(t) - a_k(s) \leq \tilde{C}_k h^\ae$, where $C, \widehat{C}_k, \tilde{C}_k$ are some non-negative constants, $\ae \in (0,1)$. In that case we can prove next theorems.

**Theorem 11.** Let stochastic process $X = \{X(t), t \in [0,T]\}$ belong to $\mathrm{Sub}_\varphi(\Omega)$ space, $\varphi(t) \in C_I$. Let $X$ be separable and assume that this process satisfies conditions (2) and **(C1)**. Then, the model $X_N(t)$ approximates the process $X(t)$ with given reliability $1 - v$ and accuracy $\delta$ in $C(0,T)$ space, if

$$v \leq 2\exp\left\{-\frac{(\zeta-1)\delta^{\frac{1}{\zeta-1}}(\delta-\gamma_N)^{\frac{\zeta}{\zeta-1}}}{\zeta(\gamma_N^\zeta(\delta-\gamma_N)+\beta^\zeta\gamma_N)^{\frac{1}{\zeta-1}}}\right\}2(e\delta)^{1/\ae},$$

when

$$\delta > \frac{\gamma_N(v^{\zeta-1}+1)+\sqrt{\gamma_N^2(v^{\zeta-1}+1)^2 + 4v^{\zeta-1}(C^\zeta(T/2)^{\ae\zeta}-\gamma_N^2)}}{2v^{\zeta-1}},$$
$$v = \min\{\beta, \gamma_N\},$$

$$\gamma_N^\zeta \leq \left(\sum_{k=1}^N \tau_\varphi^2(\xi_k)\sup_{u\in[0,T]}\delta_k^2(u) + \sum_{k=N+1}^\infty \tau_\varphi^2(\xi_k)\sup_{u\in[0,T]}a_k^2(u)\right)^{\zeta/2}.$$

P r o o f. This theorem follows from theorems 6 and 9. □

**Theorem 12.** Let stochastic process $X = \{X(t), t \in [0,T]\}$ belong to $\mathrm{Sub}_\varphi(\Omega)$ space, $\varphi(t) \in C_{II}$. Let $X$ be separable and assume that this process satisfies conditions (2) and **(C1)**. Then, the model $X_N(t)$ approximates the process $X(t)$ with given reliability $1 - v$ and accuracy $\delta$ in $C(0,T)$ space, if

$$\nu \leq 2\exp\left\{-\frac{(\zeta-1)\delta^{\frac{1}{\zeta-1}}(\delta-\gamma_N)^{\frac{\zeta}{\zeta-1}}}{\zeta(\gamma^\zeta(\delta-\gamma)+\beta^\zeta\gamma)^{\frac{1}{\zeta-1}}}\right\}2(e\delta)^{1/\text{æ}},$$

where

$$\gamma^\zeta \leq \sum_{k=1}^N \tau_\varphi^\zeta(\xi_k)\sup_{u\in[0,T]}\delta_k^\zeta(u) + \sum_{k=N+1}^\infty \tau_\varphi^\zeta(\xi_k)\sup_{u\in[0,T]}a_k^\zeta(u).$$

P r o o f. This theorem follows from the previous theorem and theorem 10. □

### Conclusions

In this paper, reliability and accuracy models of φ-sub-Gaussian stochastic processes that approximate stochastic processes in $C(T)$ spaces are considered. The case is studied when representation of a φ-sub-Gaussian stochastic process in series cannot be found explicitly and requires application of series' elements approximations. Influence of the error of such approximations is studied for reliability and accuracy of a model of stochastic process in $C(T)$ space. Theorems that allow developing a model that approximates φ-sub-Gaussian stochastic process with given reliability and accuracy in $C(T)$ space in different cases are proved.

In our further study we plan to apply the theory of generating functions to the theorems proved in current paper and study reliability and accuracy of modeling of stochastic processes in $C(T)$ space in this case.

О.М. Моклячук

ОЦІНКА НАДІЙНОСТІ ТА ТОЧНОСТІ МОДЕЛЕЙ φ-СУБ-ГАУССОВИХ ВИПАДКОВИХ ПРОЦЕСІВ У ПРОСТОРАХ $C(T)$

**Проблематика.** У теорії моделювання випадкових процесів досі не вивчались випадки побудови моделей процесів у $C(T)$, елементи розкладу яких у вигляді ряду з незалежними членами неможливо знайти у явному вигляді.

**Мета дослідження.** Метою дослідження є вивчення надійності та точності моделей процесів із просторів $\mathrm{Sub}_\varphi(\Omega)$, які не можуть бути розкладені в ряд явно, у просторі $C(T)$.

**Методика реалізації.** На базі попередніх досліджень теорії моделювання випадкових процесів із заданими надійністю та точністю розглядається припущення про можливість розкладу випадкового процесу в ряд з елементами, знайденими з певним наближенням.

**Результати дослідження.** Вивчено вплив похибки наближення елементів розкладу випадкового процесу в ряд із незалежними членами на надійність і точність побудови моделі такого процесу в просторі $C(T)$.

**Висновки.** Доведено теореми, які дадуть можливість оцінювати надійність і точність побудови моделі випадкового процесу в просторі $C(T)$ у випадку, коли елементи розкладу такого процесу в ряд із незалежними членами можуть бути знайдені лише з певною похибкою, із застосуванням, наприклад, числових методів.

**Ключові слова:** випадкові процеси; φ-суб-Гауссові процеси; моделі випадкових процесів; точність і надійність моделювання випадкових процесів.

А.М. Моклячук

ОЦЕНКА НАДЕЖНОСТИ И ТОЧНОСТИ МОДЕЛЕЙ φ-СУБ-ГАУССОВСКИХ СЛУЧАЙНЫХ ПРОЦЕССОВ В ПРОСТРАНСТВАХ $C(T)$

**Проблематика.** В теории моделирования случайных процессов на сегодняшний день не изучались случаи построения моделей процессов в $C(T)$, элементы разложения которых в виде ряда с независимыми членами невозможно найти в явном виде.

**Цель исследования.** Целью исследования является изучение надежности и точности моделей процессов из пространства $\mathrm{Sub}_\varphi(\Omega)$, которые не могут быть представлены в виде ряда явно, в пространстве $C(T)$.

**Методика реализации.** На базе предыдущих исследований в теории моделирования случайных процессов с заданными надежностью и точностью рассматривается допущение о возможности разложения случайного процесса в ряд с элементами, найденными с некоторым приближением.

**Результаты исследования.** Изучено влияние погрешности приближения элементов разложения случайного процесса в ряд с независимыми членами на надежность и точность построения модели такого процесса в пространстве $C(T)$.

**Выводы.** Доказаны теоремы, которые позволят оценивать надежность и точность построения модели случайного процесса в пространстве $C(T)$ в случае, когда элементы разложения такого процесса в ряд с независимыми членами могут быть найдены лишь с некоторой погрешностью, с применением, например, численных методов.

**Ключевые слова:** случайные процессы; φ-суб-Гауссовские процессы; модели случайных процессов; точность и надежность моделирования случайных процессов.